\documentclass[12pt]{article}
\usepackage{amsmath,amsfonts}
\usepackage[mathscr]{eucal}
\usepackage{amssymb,amsmath}
\usepackage{latexsym}
\usepackage{amsfonts}
\usepackage{amscd}
\usepackage{amsthm}
\usepackage{graphics}
\usepackage{epsfig}

\newcommand {\del}{\partial}

\newcommand {\cH}{\mathcal {H}}

\newcommand{\diam}{{\rm diam}}

\newcommand{\dist}{{\rm dist}}
\newcommand{\ball}{{\rm Ball}}

\newtheorem{theorem}{Theorem}[section]
\newtheorem{lemma}[theorem]{Lemma}

\newtheorem{cor}[theorem]{Corollary}
\newtheorem{rem}[theorem]{Remark}
\numberwithin{equation}{section}

\newlength{\originalbase}
\begin{document}
%
%

\title{Big-Pieces-of-Lipschitz-Images Implies a Sufficient Carleson Estimate in a Metric Space}

\author{Raanan Schul}

\date{}

\maketitle

\psdraft

\begin{abstract}
This note is intended to be a supplement to the bi-Lipschitz decomposition of Lipschitz maps shown in 
\cite{my-lip-bilip}.  
We show that in the case of  1-Ahlfors-regular sets, the condition of having `Big Pieces of bi-Lipschitz Images' (BPBI) is equivalent to a Carleson condition.   
\end{abstract}
\section{Introduction}
This note is intended to be a supplement to the bi-Lipschitz decomposition of Lipschitz maps shown in 
\cite{my-lip-bilip}.    
We assume familiarity with some methods from \cite{RS-metric} and some definitions from \cite {DS} (also in \cite{my-lip-bilip}).

We prove the following theorem.
\begin{theorem}\label{t:BPLI-to-GL}
Let $E$ be a 1-Ahlfors-regular set.  Suppose that $E$ has Big Pieces of Lipschitz Images.  Then for any $x\in E$ and $r<\diam(E)$ we have
\begin{equation}
\sum\limits_{B\in \widehat{\mathcal{G}}^E\atop B\subset \ball(x,2r)}
	\iiint_{(B\cap E)^3}\del(\{x_1,x_2,x_3\})\diam(B)^{-3}
								d\cH^1(x_3)d\cH^1(x_2)d\cH^1(x_1)\lesssim r\,,
\label{e:BPLI-to-GL-ver-2}
\end{equation}
where $\widehat{\mathcal{G}}^E$ is as defined as in \cite{RS-metric} and the constant $A$ in the definition of 
$\widehat{\mathcal{G}}^E$ is large enough.
\end{theorem}
\begin{rem}
As one may expect, the main point about Lipschitz images, is that we have \eqref{e:BPLI-to-GL-ver-2} for 1-Ahlfors-regular Lipschitz images.  
It is apparent  from the proof that one may replace the condition BPLI by BP(*), where (*) is any collection of 1-Ahlfors-regular sets which satisfies \eqref{e:BPLI-to-GL-ver-2}.
\end{rem}

From this, together with  what appears in \cite{my-lip-bilip} and \cite{Ha-2}, one trivially concludes  the following. 
\begin{cor}
Let $E$ be a 1-Ahlfors-regular set. 
Assume that the constant $A$ in the definition of 
$\widehat{\mathcal{G}}^E$ is large enough.
Then the following conditions are equivalent
\begin{itemize}
\item For any $x\in E$ and $r<\diam(E)$ we have the Carleson estimate \eqref{e:BPLI-to-GL-ver-2}.
\item $E$ has Big Pieces of Lipschitz Images
\item $E$ has Big Pieces of bi-Lipschitz Images
\end{itemize}
\end{cor}
This corollary is the motivation for his essay.  In general, getting a similar result to this corollary for $k-$Ahlfors-David-regular sets was a large part of the motivation for  \cite{my-lip-bilip}.  
Unfortunately, we are unable to do this for $k>1$.  
The obstacle is finding a `correct' Carleson condition, and thus far we have encountered some technical obstacles in our attempts.  
We conjecture that one can define a Carleson condition so that the above corollary will hold for $k-$Ahlfors-David-regular sets with $k>1$, thus overcoming these technical difficulties. 
In general the author is interested in pushing the David-Semmes theory of uniform rectifiability into the setting of metric spaces.  We see this theory as  related to embedability questions widely studied in the theoretical computer science community as well as related to mathematical applications  (in particular, analysis of data sets).

\begin{proof}[Proof of Theorem \ref{t:BPLI-to-GL}]
We first note that without loss of generality the we may replace  equation \eqref{e:BPLI-to-GL-ver-2} by
\begin{equation}
\sum\limits_{Q\in \Delta^i(E)\atop Q\subset Q_0}
	\ \ \iiint\limits_{(E\cap Q)\times (E\cap Q)\times (E\cap Q)}\frac{\del(\{x_1,x_2,x_3\})}{\diam(Q)^{3}}
								d\cH^1(x_3)d\cH^1(x_2)d\cH^1(x_1)\lesssim\diam(Q_0)\,,
\label{e:BPLI-to-GL-ver-3}
\end{equation}
where $i$ satisfies $1\leq i\leq P_1$, $\Delta^i(E)$ is a dyadic filtration   constructed from $\widehat{\mathcal{G}}^E$
and $Q_0$ in $ \Delta^i$ is an arbitrarily chosen `cube'.
This uses Ahlfors-regularity.
(For such constructions see   \cite{Ch,Da}.  See also  constructions in section 3.3 of \cite{RS-metric}.  Such constructions were also used in say \cite{paul-muller-book}). 
Note that $P_1$ depends only on the Ahlfors-regularity constant of $E$.

We will now follow the outline of the proof of Theorem IV.1.3  in \cite{DS}.
First we need a John-Nirenberg-Str\"ombrg  type Lemma.
This lemma is stated and proved in section IV.1.2 of \cite{DS}. (the constant coming out of the proof is $\frac{N}{\eta^2}$), and so we only state it:
\begin{lemma}\label{l:JNS}
Let $\Delta$ be a dyadic filtration.
Let $\alpha:\Delta\to [0,\infty)$ be given.  Suppose that for some $N>0$ and $\eta>0$ we have
for all $Q_0\in \Delta$
\begin{equation}
\left|  \left\{   x\in Q_0: \sum\limits_{Q\ni x \atop Q\subset Q_0, \ Q\in \Delta} \alpha(Q) \leq N \right\}  \right|
\geq \eta |Q_0|\,.
\end{equation}
Then 
\begin{equation}
\sum\limits_{ Q\subset Q_0\atop Q\in \Delta} \alpha(Q) |Q| \lesssim_{N,\eta} |Q_0|
\end{equation}
for all $Q_0\in \Delta$.
\end{lemma}

Let  $i\leq P_1$ and $Q_0\in \Delta^i$   be given.
Let $\tilde E=\tilde E (Q_0)$ be a set such 
\begin{equation}\label{e:a-big-piece}
E\cap \tilde E\cap Q_0\geq \theta \diam(Q_0)\,, 
\end{equation}
and 
$\tilde E$ is a Lipschitz Image as in the definition of Big Pieces of Lipschitz Images (BPLI).
We may assume that $\tilde E$ is actually a biLipschitz image of a subset of the real line by \cite{my-lip-bilip}.
Hence, by extending this map and the results of \cite{RS-metric}, we have (for any $A'$)
\begin{equation}\label{tilde-sum-carleson}
\sum\limits_{Q\in \Delta^i(\tilde E)\atop Q\subset Q_0}
	\ \ \iiint\limits_{\{x\in \tilde E: \dist(x,Q)\leq A'\diam(Q)\}^3}\frac{\del(\{x_1,x_2,x_3\})}{\diam(Q)^{3}}
								d\cH^1(x_3)d\cH^1(x_2)d\cH^1(x_1)
\lesssim \diam(Q_0)\,.
\end{equation}
We note that \eqref{tilde-sum-carleson} implies (by Ahlfors-regularity of $E$ and the constant $A'$  being large enough)
\begin{equation}\label{sum-carleson}
\sum\limits_{Q\in \Delta^i( E)\atop Q\subset Q_0}
	\ \ \iiint\limits_{\{x\in \tilde E: \dist(x,Q)\leq \diam(Q)\}^3}\frac{\del(\{x_1,x_2,x_3\})}{\diam(Q)^{3}}
								d\cH^1(x_3)d\cH^1(x_2)d\cH^1(x_1)
\lesssim \diam(Q_0)\,.
\end{equation}
We  have for any $Q_1$ in $\Delta_i$
\begin{equation}
\begin{aligned}
\sum\limits_{Q\subset Q_1,\ Q\in \Delta^i( E) \atop 
Q\cap\tilde E\neq \emptyset} 
\ 	\int\limits_{Q\cap E} \frac{\dist(\cdot, \tilde E)}{\diam(Q)} 
&\leq&
\int\limits_{x\in Q_1\cap E}  \sum\limits_{Q\ni x,\ Q\in \Delta^i( E) \atop  Q\in Q_1,\ Q\cap \tilde E \neq \emptyset} 
		\frac{\dist(\cdot, \tilde E)}{\diam(Q)}\\
&\lesssim&
\int\limits_{x\in Q_1\cap E} 1+\frac12 +\frac14 +\frac18 + ...
\lesssim
\diam(Q_1)\,.
\end{aligned}
\end{equation}
Similarly to the proof of Lemma 3.11 in \cite{RS-metric}, we have
\begin{equation}
\begin{aligned}
&\ \ \iiint\limits_{(E\cap Q)\times (E\cap Q)\times (E\cap Q)}\frac{\del(\{x_1,x_2,x_3\})}{\diam(Q)^{3}}
								d\cH^1(x_3)d\cH^1(x_2)d\cH^1(x_1)\\		
\lesssim&
\ \ \int\limits_{E\cap Q} \frac{\dist(\cdot, \tilde E)}{\diam(Q)}\ +\\
&
\iiint\limits_{\{x\in\tilde E: \dist(x,Q)\leq \diam(Q)\}^3}\frac{\del(\{x_1,x_2,x_3\})}{\diam(Q)^{3}}
								d\cH^1(x_3)d\cH^1(x_2)d\cH^1(x_1)\,.
\end{aligned}\label{e:beta-approx}
\end{equation}
Summing equation \eqref{e:beta-approx} we get
\begin{equation}
\sum\limits_{Q\subset Q_0,\ Q\in \Delta^i( E) \atop 
Q\cap\tilde E\neq \emptyset} 
		\ \ \iiint\limits_{(E\cap Q)\times (E\cap Q)\times (E\cap Q)}
		\frac{\del(\{x_1,x_2,x_3\})}{\diam(Q)^{3}}
								d\cH^1(x_3)d\cH^1(x_2)d\cH^1(x_1)
\leq C \diam(Q_0)		
\end{equation}
From \eqref{e:a-big-piece} and taking 
$N=2C$ and $\eta=\frac12\theta$ and 
$$\alpha(Q)=\ \ \iiint\limits_{\{x\in E: \dist(x,Q)\leq \diam(Q)\}^3}
	\frac{\del(\{x_1,x_2,x_3\})}{\diam(Q)^{4}}
								d\cH^1(x_3)d\cH^1(x_2)d\cH^1(x_1)$$
we have the conditions for Lemma \ref{l:JNS}, which gives us Theorerm \ref{t:BPLI-to-GL}.
\end{proof}

\bibliographystyle{alpha}
\bibliography{../../bibliography/bib-file}   

\begin{thebibliography}{Dav91}

\bibitem[Chr90]{Ch}
Michael Christ.
\newblock A {$T(b)$} theorem with remarks on analytic capacity and the {C}auchy
  integral.
\newblock {\em Colloq. Math.}, 60/61(2):601--628, 1990.

\bibitem[Dav91]{Da}
Guy David.
\newblock {\em Wavelets and singular integrals on curves and surfaces}, volume
  1465 of {\em Lecture Notes in Mathematics}.
\newblock Springer-Verlag, Berlin, 1991.

\bibitem[DS93]{DS}
Guy David and Stephen Semmes.
\newblock {\em Analysis of and on uniformly rectifiable sets}, volume~38 of
  {\em Mathematical Surveys and Monographs}.
\newblock American Mathematical Society, Providence, RI, 1993.

\bibitem[Hahar]{Ha-2}
Immo Hahlomaa.
\newblock Curvature and {L}ipschitz parametrizations in 1-regular metric
  spaces.
\newblock {\em Annales Academiae Scientiarum Fennicae}, To appear.

\bibitem[M{\"u}l05]{paul-muller-book}
Paul F.~X. M{\"u}ller.
\newblock {\em Isomorphisms between {$H\sp 1$} spaces}, volume~66 of {\em
  Instytut Matematyczny Polskiej Akademii Nauk. Monografie Matematyczne (New
  Series) [Mathematics Institute of the Polish Academy of Sciences.
  Mathematical Monographs (New Series)]}.
\newblock Birkh\"auser Verlag, Basel, 2005.

\bibitem[Sch]{my-lip-bilip}
Raanan Schul.
\newblock Bi-{L}ipschitz decomposition of {L}ipschitz functions into a metric
  space.
\newblock {\em {\tt arXiv:math/0702630}}.

\bibitem[Schar]{RS-metric}
Raanan Schul.
\newblock Ahlfors-regular curves in metric spaces.
\newblock {\em Annales Academiae Scientiarum Fennicae}, To appear.

\end{thebibliography}

\end{document}